\documentclass[english]{amsart}
\usepackage[T1]{fontenc}
\usepackage[latin9]{inputenc}

\usepackage{babel}
\usepackage{amsmath, amssymb}
\newtheorem{Theorem}{Theorem}[section]

\newtheorem{Lemma}[Theorem]{Lemma}
\newtheorem{Proposition}[Theorem]{Proposition}
\newtheorem{Definition}[Theorem]{Definition}
\newtheorem{rem}[Theorem]{Remark}
\newtheorem{ex}[Theorem]{Example}
\numberwithin{equation}{section}

\newcommand{\R}{\mathbb R}
\newcommand{\N}{\mathbb N}

\newcommand{\eps}{\varepsilon}

\newcommand{\C}{\mathcal{C}}
\newcommand{\F}{\mathcal{F}}

\newcommand{\X}{\mathcal{X}}
\begin{document}

\title{The mean value for infinite volume measures, infinite products and heuristic infinite dimensional Lebesgue measures}

\author{Jean-Pierre Magnot}
\address{Lycée Jeanne d'Arc; Avenue de Gande Bretagne; F-63000 Clermont-Ferrand}
\email{jean-pierr.magnot@ac-clermont.fr}
\begin{abstract}
one of the goals of this article is to define a an unified setting adapted to the description of means (normalized integrals or invariant means) on an infinite product of measured spaces with infinite measure.
We first remark that some known examples coming from the theory of metric measured spaces and also from oscillatory integrals are obtained as limits 
of means with respect to finite measures. Then, we explore in a systematic
way the limit of means of the type 
$$ \lim \frac{1}{\mu(U_n)} \int_{U_n} f d\mu$$
where $\mu$ is a 
 a $\sigma-$finite Radon measure $\mu.$ 
In some cases, we get a linear extension of the limit at infinity. 
Then, the mean value on an infinite product is defined, 
first for cylindrical functions and secondly taking the uniform limit. 
Finally, the mean value for the heuristic Lebesgue measure 
on a separable  infinite dimensional topological vector space 
(but principally on a Hilbert space) is defined. This last object is shown to be invariant by translation, scaling and restriction.

\end{abstract}
\maketitle

\tableofcontents

\section*{Introduction}

The very early starting point of this work is the well-known lack of
adequate definition of an  
infinite dimensional Lebesgue measure on a Hilbert space. 
Such a measure, which is assumed invariant by translation and 
by action of the orthogonal or unitary group, 
is often said non existing in the infinite dimensional setting, because
it has not enough sets with finite measure, see e.g. \cite{Ba1,Ba2} for an overview.  
Anyway, there exists many measures on infinite dimensional objects, but none 
can be used to replace the Lebesgue measure in mathematical constructions. 
From another point of view, normalized infinite dimensional integrals are
well defined objects, and in some sense appear as limits of finite measures, see e.g. \cite{AHM}.
Yet in another setting, there is a property of concentration of measure in metric measured spaces in which can coincide with the definition of a mean for uniform functions in e.g. $\mathcal{S}^\infty$, see e.g. \cite{Gro,GroMi,Pes}.
These two last approaches appear as relevant of the same procedure: defining means from limits of measures. This is why, including a previous work \cite{MaICM}, we suggest a general setting in section \ref{gen} for means spanned by finite measures, where a measure is just a tool to define
the \textbf{mean value} of the function, namely
$$\frac{1}{\mu(X)}\int f d\mu$$
for a set $X$ of $\mu-$finite measure.

\vskip 12pt

The theory developed in section \ref{gen} is then specialized
to a restricted class of means, first to the means obtained with a $\sigma-$finite Radon measure, using a creasing sequence $(U_n)_{n\in\N}$
 of Borel subsets with finite measure satisfying 
$\bigcup_{n \in \N}U_n = X$ among other technical conditions, by:
$$ \bar{f} = \lim_{n \rightarrow +\infty}  \frac{1}{\mu(U_n)} \int_{U_n} f d\mu,$$
in the spirit of convergence of finite dimensional oscillatory integrals.  
In section
\ref{MS}, we develop the basics of this theory on a measured space. Since this mean value depends (in general)
on the sequence $U$ and on the measure $\mu,$ we do not adopt the notation $\bar{f}$ but prefer $WMV_\mu^U(f)$ or $MV_\mu(f)$, 
abbreviations for ``Weak Mean Value'' and for `` Mean Value''. Formulas for changing of measure leads us to 
an extension of the asymptotic comparison of functions ($f\sim_\infty g,$ $f = O(g)$ and $f=o(g)$ ) to measures. 
As a particular case, the mean value with respect to the Lebesgue measure on $\R$ appears as a linear extension 
of the limit at $\infty$ of functions. 
We know very few about the behaviour of the mean value of limit of functions: 
the mean value is not continuous for vague convergence, 
but continuous for uniform convergence. There is certainly an intermediate kind of 
convergence more adapted to mean values, to be determined. 
 We also give an application of this notion: the homology map
as a mean value of a function on the space of harmonic forms, using Hodge theory.  

Secondly, we get to infinite products of measured spaces in section \ref{IP}. Recall that there is an induced measure 
on an infinite product of measured spaces only if we have spaces with finite measures. 
Our approach here is mostly inspired by Daniell's integral, 
which is a preliminary approach to Wiener measure. We consider cylindrical functions, 
and define very easily their mean values as mean values of functions defined on a finite product 
of measured space. Then, we extend it to functions that are uniform limits of sequences of cylindrical functions. 
As an application, we give a definition of the mean value on infinite configuration spaces for Poisson measure. 

Finally, we get to vector subspaces of Hilbert spaces in section 3. 
This is where we decide to focus on the announced heuristic 
infinite dimensional Lebesgue measure. 
The mean value is developed and we study its invariance properties. 
It appears invariant by translation and by scaling, 
and also by action of the unitary group. 
But the last one remains dependent on the choice of the orthonormal 
basis used for the definition, which is analogous to the multiplicative 
anomaly of renormalized determinants (see e.g.\cite{KV} 
for the canonical determinant of Kontsevich and Vishik) since it can be read as 
a non invariance while changing the basis.
 As a concluding remark, we show that 
this approach has a technical difference with the approach by measures on infinite dimensional spaces. 
We show that the mean value of a bounded continuous function $f$ remains the same while 
restricting to a dense vector subspace. This exhibits a striking difference 
from e.g. the Wiener measure on continuous paths, for which the space of $H^1$ paths is of measure $0.$  
With all these elements, we can now explain where is the originality of our approach. Here, the total 
volume is not considered as a constant of the total space, 
but as a scale-like element 
to compare with the integral of a function. 
This is exactly the spirit of the formula of the mean value in finite volume. 
 
\section{The space of means spanned by sequences of finite measures}\label{gen}

Let $(X,\mu)$ be a measured space. Following \cite{Pa}, \cite{Pes}, let us fix a vector subspace $\mathcal{F}\subset L^\infty(X,\mu)$ such that $1_X \in \mathcal{F}.$ A \textbf{mean} on $\mathcal{F}$ is a linear map $\phi:\mathcal{F}\rightarrow \mathbb{C}$
such that $\phi(1_X)=1.$ Alternately, if $(X,d)$ is a metric space, 
given $\mathcal{F}\subset C^0_b(X)$ (space of bounded maps), 
a \textbf{mean} on $\mathcal{F}$ is a linear map 
$\phi:\mathcal{F}\rightarrow \mathbb{C}$ such that $\phi(1_X)=1.$
These two terminologies come from the basic example where $\mu$ is a Borel
probability measure on a metric space $(X,d),$ for which the mean of a continuous integrable map $f$ is its expectation value $$\int_X f d\mu,$$
and can be approximated by sequences of barycenters of Dirac measures via Monte Carlo methods. We intent to decribe means spanned by approximations via finite measures on a metric space in this section.

\subsection{Means spanned by probability measures}

Let $X$ be a complete metric space and let $C^0_b(X)$ be the space of bounded complex valued continuous maps on $X.$ We note by $\mathbb{P}(X)$ the space of Borel probability measures on $X.$ 
\begin{Definition}
A complex (resp. real) \textbf{probability mean} is a linear map $\tau:\mathcal{D}_\tau \subset C^0_b(X)\rightarrow \C $ which is defined as the limit of barycenters with complex  (resp. real) weights of a sequence of Borel probability measures on $X,$ i.e.
for $\mathbb{K}=\mathbb{R}$ or $\mathbb{C},$
$$\exists (\mu_n, \alpha_n)_{n \in \N}\in (\mathbb{P}(X) \times \mathbb{K})^\N, \forall m\in \N^*, $$ $$\left\{ \sum_{n = 0}^{m} \alpha_n \neq 0 \right\} \quad \wedge \quad \left\{ \forall f \in C^0_b(X), \tau(f) = \lim_{m \rightarrow +\infty} \frac{1}{\sum_{n = 0}^{m} \alpha_n} \left( \sum_{n = 0}^{m} \alpha_n {\mu_n}(f) \right)\right\}.$$
\end{Definition}

We note by $\widetilde{\mathcal{PM}}_\mathbb{K}(X)$ the space of $\mathbb{K}-$probability means, by ${\mathcal{PM}}_\mathbb{K}(X)$ the set of probability means $\tau$ such that $\mathcal{D}_\tau = C^0_b(X),$ by $\widetilde{\mathcal{PM}}_\mathbb{R}^+(X)$ the means $\tau$ obtained by a sequence $(\alpha_n)_{n \in \N}\in \R_+^*$ and by ${\mathcal{PM}}_\mathbb{R}^+(X)$ the space $\widetilde{\mathcal{PM}}_\R(X)\cap\widetilde{\mathcal{PM}}_\mathbb{R}^+(X) $  We have a spacial class spanned by the Dirac measures:

\begin{Definition} \cite{MaICM}
A complex (resp. real) \textbf{Dirac mean} is a linear map $\tau:\mathcal{D}_\tau \subset C^0_b(X)\rightarrow \C $ which is defined as the limit of barycenters with complex  (resp. real) weights of a sequence of Dirac measures on $X,$ i.e.
for $\mathbb{K}=\mathbb{R}$ or $\mathbb{C},$
$$\exists (x_n, \alpha_n)_{n \in \N}\in (X \times \mathbb{K})^\N, \forall m\in \N^*, $$ $$\left\{ \sum_{n = 0}^{m} \alpha_n \neq 0 \right\} \quad \wedge \quad \left\{ \forall f \in C^0_b(X), \tau(f) = \lim_{m \rightarrow +\infty} \frac{1}{\sum_{n = 0}^{m} \alpha_n} \left( \sum_{n = 0}^{m} \alpha_n \delta_{x_n}(f) \right)\right\}.$$
\end{Definition}
 
We note by $\widetilde{\mathcal{DM}}_\mathbb{K}(X),$  ${\mathcal{DM}}_\mathbb{K}(X), $   $\widetilde{\mathcal{DM}}_\mathbb{R}^+(X),$  ${\mathcal{DM}}_\mathbb{R}^+(X) $ the sets of Dirac means corresponding respectively to  $\widetilde{\mathcal{PM}}_\mathbb{K}(X),$  ${\mathcal{PM}}_\mathbb{K}(X), $   $\widetilde{\mathcal{PM}}_\mathbb{R}^+(X),$  ${\mathcal{PM}}_\mathbb{R}^+(X) $

\begin{Proposition}
$\widetilde{\mathcal{PM}}_\mathbb{K}(X),$ $\mathcal{PM}_\mathbb{K}(X)$
$\widetilde{\mathcal{DM}}_\mathbb{K}(X)$ and $\mathcal{DM}_\mathbb{K}(X)$ are $\mathbb{K}-$affine spaces.
\end{Proposition}
The proof is obvious adapting elementary proofs on (classical, finite) barycenters. We give also the following, in order to make quickly the link
with the Monte-Carlo method.

\begin{Proposition}
If $X$ is moreover a locally compact manifold, one has the following inclusions:
 
\begin{itemize}
\item $\mathbb{P}(X) \subset \mathcal{DM}_\mathbb{R}^+(X).$
\item If $X$ is compact, $\mathbb{P}(X) = \mathcal{DM}_\mathbb{R}^+(X)= \mathcal{PM}_\mathbb{R}^+(X).$
\end{itemize}
\end{Proposition}

\noindent\textbf{Proof.}
\begin{itemize}
\item Let $\mu \in \mathbb{P}(X)$ and let $(x_n)_{n\in  \N}$
be a uniformly distributed sequence with respect to $\mu.$
Then, $\forall f\in C^0_b(X),   \lim_{m\rightarrow +\infty} \frac{1}{n+1} \sum_{n=0}^m f(x_n) = \mathbb{E}_\mu(f) = \int_X f d\mu.$
Thus
$$\mu(f) = \lim_{m\rightarrow +\infty} \frac{1}{n+1} \sum_{n=0}^m \delta_{x_n}(f).$$
\item If $X$ is compact, the space of (signed) finite measures on $X$ coincide with the
$\left(C^0_b(X)\right)'.$ Since $\mathcal{PM}_\mathbb{R}(X)\subset\left(C^0_b(X)\right)',$
we easily get the result.
\end{itemize} \qed

\subsection{Limit means}
\begin{Definition}\label{lmean}
Let $X=(X_n,\tau_n)_{n \in \N}$ be a sequence of probability spaces such that 

- $\forall n \in \N, $ $X_n$ is a metric space.

- $\forall n \in \N, X_n \subset X_{n+1},$ and the topology of $X_{n+1}$ restricted to $X_n$ co\"incides with the topology of $X_n.$

- $\forall n \in \N, \tau_n\in \widetilde{\mathcal{PM}}_\C (X_n).$
Then, we define, for the maps $f$ defined on $\bigcup_{n \in \N} X_n$ with values in a complete topological vector space, if $\forall n \in \N, f_{|X_n} \in \mathcal{D}_\tau$ and if the limit converges, 
$$ LM^X(f) =\lim_{n\rightarrow +\infty} \tau_n (f)$$
called
\textbf{limit mean} of $f$ with respect to $X.$  
 
\end{Definition}

\subsection{Probability means in the mm-space setting}
We use two handbooks for preliminaries on these notions: \cite{Gro}
and \cite{Pes}. 
\begin{Definition}\cite{GroMi}
A \textbf{space with metric and measure}, or a \textbf{metric measured space}
(mm-space for short) is a triple $(X, d, \mu)$ where $(X,d)$
is a metric space and $\mu$ is a probability measure on the
Borel tribu on $X.$ 
\end{Definition}

Let $A \subset X,$ let $\varepsilon >0.$ We note by
$$A_\eps = \{ x \in X | d(A,x)<\eps \}.$$

\begin{Definition}\cite{GroMi}
A \textbf{L\'evy family} is a sequence 
$\mathcal{X}=(X_n, d_n, \mu_n)_{n \in \N}$ of mm-spaces
if, for each sequence $(A_n)_{n \in \N}$ such that 
$$\forall n \in \N, A_n \hbox{ is a Borel subset of } X_n$$
ans satisfying 
$$ \hbox{lim inf}_{n\rightarrow +\infty} \mu_n(A_n) > 0,$$
then 
$$\forall \eps > 0, \lim_{n \rightarrow +\infty}
\mu_n\left((A_n)_\eps\right)= 1.$$ 
\end{Definition}

\vskip 12pt

In the sequel, we shall assume that 
$$ \forall n \in \N, X_n \subset X_{n+1}$$
with continuous injection. Notice that we do not assume that
$d_n$ is the restriction of $d_{n+1}$ hich allows us some 
freedom on metric requirements. The technical necessary condition is the following: let $n \in \N$ and let $B_{n+1}$ be a Borel subset of $X_{n+1}.$ Then $B_{n+1} \cap X_n$ is a Borel subset of $X_n.$ 
We have here a priori a class of limit means following the terminology of Definition \ref{lmean}.
 Let us quote first the classical (and historical) example of a Levy familysee e.g. \cite{Gro}, section 3$\frac{1}{2}$.19,
which gives an example of mean value:

\begin{ex} {\bf The Levy family of spheres and the concentration phenomenon}

Let us consider the seuquence of inclusions
$$ S^1 \subset S^2 \subset ... \subset S^n \subset S^{n+1} ... \subset S^\infty=\bigcup_{n = 1}^\infty S^n$$
equipped with the classical Euclide (or Hilbert) distance and (except for $S^\infty$) the normalized spherical measure $\mu$ (we drop the index for the measure in sake of clear notations). Then, for any $\R -$valued 1-Lipschitz function on $S^\infty,$ there exists $a\in \R$ such that:
$$\forall \epsilon >0, \quad \mu\left\{ x \in S^n | ||f(x) - a||>\epsilon \right\} < 2e^{-\frac{(n-1)\epsilon^2}{2}}.$$
\end{ex}
In a more intuitive formuation, one can say that any 1-Lipschitz function concentrates around a real vaule $a$ with respect to $\mu.$ We leave the reader with the reference \cite{Gro} for more on the metric geometry of this example.
We can reformulate: 
\begin{Proposition}
Let $\mathcal{X}=(S^n; ||.||; \mu)_{n \in \N^*}$. Then 
for any 1-Lipschitz function $f$ defined on $S^\infty,$
and with the notations used before,
$$LM^\mathcal{X}(f)=a.$$
\end{Proposition} 
\vskip 12pt
\noindent
\begin{ex}{ L\'evy families induced by Lebesgue measures}
Let $m, n\in (\N^*)^2.$ 
Take $K_m \subset`\mathbb{R^n}.$ 
For each $m `\in \N^*,$ we equip $K_m$ with the usual 
distance $d$ induced by $\R^n$ and with the probability measure
$$\mu_n = \frac{\mathbf{1}_{K_n}}{\lambda(K_n)} \lambda.$$
Setting $\mathcal{K}= (K_m,d, \mu_m)_{m \in \N^*},$
we get that $\mathcal{K}$ is a L\'evy family, but there is no concentration property. This example will be studied in the next sections of this article.
\end{ex}

\begin{Definition}
Let $f: X \rightarrow \mathbb{C}$ be a map
such that for each $ n \in \N,$ the restriction of $f$ to 
$X_n$ is $\mu_n-$integrable. Then, the \textbf{mean value} 
of $f$ with respect to the family $\mathcal{X}$ is
$$WMV^\X(f)=\lim_{n \rightarrow +\infty} \int_{X_n} f d\mu_n$$
if the limit exists. 
\end{Definition}

\subsection{Means defined by oscillatory integrals}

Let $\Phi \in C^\infty(\R^n, \R)$ be a fixed function. Following \cite{ET} (see e.g. \cite{AHM}, \cite{AMaz1}, \cite{Dui}, \cite{Ste}), we define:
\begin{Definition}
Let $f$ be a measurable function on $\R^n.$
Let $\varphi \in \mathcal{S}(\R^n)$ be a weight function such that $\varphi(0)=1.$
if the limit 
$$\lim_{\epsilon \rightarrow 0}\int_{\R^n}e^{i\Phi(x)} f(x) \varphi(\epsilon x) dx$$
exists and is independent of the fixed function $\varphi,$
then this limit is called oscillatory integral of $f$ with respect to $\Phi,$
noted $$\int_{\R^n}^o e^{i\Phi(x)} f(x) dx .$$
\end{Definition} 

The choice $\Phi(x) = \frac{i}{2h} |x|^2$ is of particular interest, and is known under the name of \textbf{Fresnel integral}. This choice gives us a mean, up to normalization by a factor $(2i\pi h)^{-\frac{d}{2}},$
and can be generalized to a Hilbert space $\mathcal{H}$ the following way:

\begin{Definition}
A Borel measurable function $f: \mathcal{H}\rightarrow \C$ is called $h-$ integrable in the sense of Fresnel is for each creasing 
sequence of projectors $(P_{n})_{n \in \N}$ such that $\lim_{n \rightarrow +\infty} P_n = Id_\mathcal{H},$ the finite dimensional approximations of the oscillatory integrals of $f$ 
$$\left\{\int_{Im P_n}^o e^{\frac{i}{2h}|P_n(x)|^2} f(P_n(x)) d(P_n(x))\right\}
\left\{\int_{Im P_n}^o e^{\frac{i}{2h}|P_n(x)|^2}  d(P_n(x))\right\}^{-1}$$
are well-defined and the limit as $n\rightarrow +\infty$ does not depend on the sequence $(P_{n})_{n \in \N}.$ In this case, it is called infinite dimensional Fresnel integral of $f$ and noted 
$$\int_{\mathcal{H}}^o e^{\frac{i}{2h}|x|^2} f(x) d(x).$$
\end{Definition}

The invariance under the choices of the map $\varphi$ and the projections $P_n$ is assumed mostly to enable stronger analysis on these objects, which intend to be useful to describe physical quantities and hence can be manipulated by physicists who sometimes work ``with no fear on the mathematical rigour'' of their calculations. But we can also remark that:
\begin{itemize}
\item the map $$f \mapsto \int_{\R^n}^o e^{i\Phi(x)} f(x) dx \in \widetilde{\mathcal{PM}}_\mathbb{\mathbb{C}}(X),$$
\item the $$f \mapsto \int_{\R^n}^o e^{i\Phi(x)} f(x) dx $$ is a limit mean through the sequence $\R \subset...\subset \R^n \subset \R^{n+1} \subset ... \subset \mathcal{H}.$
\end{itemize}
The limit mean obtained is got through the classical trick of cylindrical functions, which we shall also use in the sequel.
\section{Mean value on a measured space} \label{MS}

\subsection{Definitions}

Let $(X,\mu)$be a topological space equipped with a measure $\mu.$
Let $\mathcal{T}(X)$ be the tribu on $X$. We note by $Ren_{\mu}(X)$
the set of sequences $U=(U_{n})_{n\in\mathbb{N}}\in\mathcal{T}(X)^{\mathbb{N}}$such
that

\begin{enumerate}
\item $\bigcup_{n\in\mathbb{N}}U_{n}=X$
\item $\forall n\in\mathbb{N},$$0<\mu(U_{n})<+\infty$ and $U_{n}\subset U_{n+1}.$
\end{enumerate}
\emph{Remark:}We have in particular $lim_{n\rightarrow+\infty}\mu(U_{n})=\mu(X).$

In what follows we assume the natural condition $Ren_{\mu}\neq\emptyset.$

\begin{Definition}\label{def1}

Let $U\in Ren_{\mu}.$ Let $V$ be a separable complete locally convex topological
vector space (sclctvs). Let $f:X\rightarrow V$ be a measurable map.
We define, if the limit exists, the weak mean value of $f$ with respect
to $U$ as: \[
WMV_{\mu}^{U}(f)=\lim_{n\rightarrow+\infty}\frac{1}{\mu(U_{n})}\int_{U_{n}}fd\mu\]
 Moreover, if $WMV_{\mu}^{U}(f)$does not depend on $U$, we call
it mean value of $f$, noted $MV_{\mu}(f).$

\end{Definition}

\begin{rem}
There is a well-known integration theory for measurable Banach-valued maps. 
A separable complete locally convex topological vector space can be 
seen topologically as the projective limit of a sequence of Banach spaces. 
So that, integrating a function with image in a sclctvs is just 
considering integration on Banach spaces, and after taking the projective limit.
\end{rem}
Notice that

\begin{itemize}
\item if $V=\mathbb{R},$ setting $f_{+}=\frac{1}{2}(f+|f|)$ and $f_{-}=\frac{1}{2}(f-|f|),$
$WMV_{\mu}^{U}(f)=WMV_{\mu}^{U}(f_{+})+WMV_{\mu}^{U}(f_{-})$ for
each $U\in Ren_{\mu},$ if $f,$ $f_+$ and $f_-$ have a finite mean value.
\item The same way if $V=\mathbb{C},$$WMV_{\mu}^{U}(f)=WMV_{\mu}^{U}(\Re f)+iWMV_{\mu}^{U}(\Im f)$
for each $U\in Ren_{\mu}.$
\item We note by $\mathcal{F}_{\mu}^{U}$ the set of functions $f$ such
that $WMV_{\mu}^{U}(f)$ exists in $V,$ and by $\mathcal{F}_{\mu}$the set
of functions $f$ such that $MV_{\mu}(f)$ is well-defined.
\end{itemize}
\textbf{Examples.}

\begin{enumerate}
\item Let $(X,\mu)$ be an arbitrary measured space. Let $f=1_{X}.$ Let
$U\in Ren_{\mu}.$ $\forall n\in\mathbb{N},$ $\frac{1}{\mu(U_{n})}\int_{U_{n}}fd\mu=\frac{\mu(U_{n})}{\mu(U_{n})}=1.$
So that \[
MV_{\mu}(1_{X})=1.\]

\item Let $(X,\delta_{x})$ be a space $X$ equipped with the Dirac measure
at $x\in X$. Let $f$ be an arbitrary map to an arbirary clcvs. $U\in Ren_{\delta_{x}}\Leftrightarrow\forall n\in\mathbb{N},\quad\delta_{x}(U_{n})>0\Leftrightarrow\forall n\in\mathbb{N},\quad x\in U_{n}.$Thus,
if $U\in Ren_{\delta_{x}}$ $\forall n\in\mathbb{N},$ $\frac{1}{\delta_{x}(U_{n})}\int_{U_{n}}fd\delta_{x}=f(x).$
So that $$
MV_{\delta_{x}}(f)=f(x).$$

\item Let $(X,\mu)$ be a measured space with $\mu(X)<+\infty$. Let $f$
be an arbitrary bounded measurable map. Then one can show very easily
that we recover the classical mean value of $f:$ \[
MV_{\mu}(f)=\frac{1}{\mu(X)}\int_{X}fd\mu.\]

\item Let $X=\mathbb{R}$ equipped with the classical Lebesgue measure $\lambda.$
Let $g\in L^{1}(\mathbb{R},\mathbb{R}_{+})$ (integrable $\mathbb{R}_{+}$-valued
function). Let $U\in Ren_{\lambda}.$ We have that $\lim_{n\rightarrow+\infty}\int_{U_{n}}gd\lambda\leq\int_{\mathbb{R}}gd\lambda<+\infty$
so that \[
MV_{\lambda}(g)=0.\]

\item Let $X=\mathbb{R}$ equipped with the Lebesgue measure $\lambda.$
Let $f(x)=\sin(x)$ and let $U_{n}=[-(n+1);(n+1)].$ The map $\sin$
is odd so that $WMV_{\lambda}^{U}(\sin)=0.$ Now, let $U_{n}'=[-2\pi n;2\pi n]\cup\bigcup_{j=0}^{n}[2(n+j)\pi;(2(n+j)+1)\pi].$
Then $WMV_{\lambda}^{U'}(\sin)=\frac{1}{5\pi}.$ This shows that $\sin$ has no
(strong) mean value for the Lebesgue measure.

\item Let $X=\N$ equipped with $\gamma$ the counting measure. 
Let $n \in \N$ and set $U_n = [0;n] \cap \N.$ Let $(u_n) \in \R^\N$ and 
$U= (U_n)_{n \in \N}.$ Then, 
$$WMV^U_\gamma(u_n)=\lim_{n \rightarrow +\infty} \frac{1}{n+1}\sum_{k=0}^n u_k$$
is the Cesar\`o limit. 

\end{enumerate}

\subsection{Basic properties}

In what follows and till the end of this paper we assume the natural
condition $Ren_{\mu}\neq\emptyset.$

\begin{Proposition} \label{linear}

Let $(X,\mu)$be a measured space. Let $U\in Ren_{\mu}.$Then

\begin{enumerate}
\item $\mathcal{F}_{\mu}^{U}$ is a vector space and $WMV_{\mu}^{U}$ is
linear
\item $\mathcal{F}_{\mu}$ is a vector space and $MV_{\mu}$ is linear.
\end{enumerate}
\end{Proposition}

The proof is obvious. 

We now clarify the preliminaries that are necessary to study the perturbations of the mean value of a fixed function
with respect to perturbations of the measure. 
\begin{Proposition}\label{addition}
Let $\mu$ and $\nu$ be Radon measures. Let $U=(U_n)_{n\in\N}\in Ren_\mu \cap Ren_\nu .$
Assume that $\Theta (\mu ,\nu) = \lim_{n\rightarrow + \infty} \frac{\mu (U_n) }{ (\mu + \nu)(U_n)} \in [0;1]$ exists.  
Then  $f \in\F_{\mu + \nu}^U $ and $$WMV_{\mu+\nu}^U(f) = \Theta(\mu,\nu) WMV_\mu^U(f) + \Theta(\nu, \mu)WMF_\nu^U(f).$$

\end{Proposition}

\noindent
\textbf{Proof.}
\begin{eqnarray*} \frac{1 }{(\mu + \nu)(U_n)} \int_{U_n} f d(\mu + \nu) &=& \frac{\mu(U_n) }{ (\mu + \nu)(U_n)} \left\{ \frac{1}{ \mu (U_n)} \int_{U_n} f d\mu\right\}\\ & + &\frac{\nu(U_n) }{ (\mu + \nu)(U_n)} \left\{ \frac{1 }{ \nu (U_n)} \int_{U_n} f d\nu\right\}.\end{eqnarray*}
Thus, we get the result taking the limit. 

\begin{Proposition} \label{scaling}
Let $\mu$ be a Radon measure and let $k \in \R_+^*.$ Then $Ren_{k \mu} = Ren_\mu, $ $\F_{k \mu} = \F_\mu ; $
moreover $\forall U \in Ren_\mu,$ $\F_{k \mu}^U = \F_\mu^U$ and $WMV_{k \mu}^U=WMF_\mu^U.$
\end{Proposition} 

The proof is obvious.

\begin{Theorem}
Let $\mu$ be a measure on $X,$ let $U \in Ren_\mu$ and $f\in \F_\mu^U.$ Let 
$$\mathcal{M}(\mu, U , f) = 
\left\{ \nu | U \in  Ren_\nu \hbox{ and } WMV_\nu^U(f)=WMV_\mu^U(f) \right\} $$
 
$\mathcal{M}(\mu,U,F)$ is a convex cone.

\end{Theorem}

\vskip 12pt
\noindent
\textbf{Proof.}
Let $k>0$ and let $\nu \in \mathcal{M}(\mu,U,F)$. Setting $\nu' = k\nu, $
we get $WMV^U_{\nu'}(f)=WMV^U_{\nu}(f)$ by Proposition \ref{scaling}, thus $\mathcal{M}(\mu,U,F)$ is a cone. 

Now, let $(\nu, \nu')\in \mathcal{M}(\mu,U,F)^2.$ Let $t\in [0;1]$ and let $\nu'' = t\nu + (1-t)\nu'.$ 

$\bullet$ \underline{Let us show that $U \in Ren_{\nu''}$}.  

Let $n \in \N.$ We have $\nu'' (U_n) = t\nu(U_n) + (1-t)\nu'(U_n)$, so that $\nu''(U_n) \in \R_+^*.$

$\bullet$ \underline{Let us show that $WMV_{\nu''}^U(f) = WMV_\mu^U(f)$}. 
We already know that $WMV_{\nu'}^U(f)=WMV_{\nu}^U(f)=WMV_{\mu}^U(f).$ Let $n \in \N.$
\begin{eqnarray*} \frac{1 }{\nu''(U_n)} \int_{U_n} f d(\nu'') & = & 
\frac{t\nu(U_n) }{ (t\nu + (1-t)\nu')(U_n)} \left\{ \frac{1}{ t\nu (U_n)} \int_{U_n} f d(t\nu)\right\} \\
 &&+ \frac{(1-t)\nu'(U_n) }{ (t\nu + (1-t)\nu')(U_n)} \left\{ \frac{1 }{ (1-t)\nu' (U_n)} 
\int_{U_n} f d((1-t)\nu')\right\} \\
& = & \left\{\frac{t\nu(U_n) }{ (t\nu + (1-t)\nu')(U_n)} + \frac{(1-t)\nu'(U_n) }{ (t\nu + (1-t)\nu')(U_n)} \right\} WMV^U_\mu(f)\\
&&+
\frac{t\nu(U_n) }{ (t\nu + (1-t)\nu')(U_n)} \left\{ \frac{1}{ \nu (U_n)} \int_{U_n} f d(\nu)-WMV^U_\mu(f)\right\} \\
 &&+ \frac{(1-t)\nu'(U_n) }{ (t\nu + (1-t)\nu')(U_n)} \left\{ \frac{1 }{ \nu' (U_n)} 
\int_{U_n} f d\nu'-WMV^U_\mu(f)\right\}
\end{eqnarray*}
Now, we remark that 
$$ \left\{\frac{t\nu(U_n) }{ (t\nu + (1-t)\nu')(U_n)} + \frac{(1-t)\nu'(U_n) }{ (t\nu + (1-t)\nu')(U_n)} \right\}=1,$$
and that 
$$\lim_{n \rightarrow + \infty} \frac{t\nu(U_n) }{ (t\nu + (1-t)\nu')(U_n)} \left\{ \frac{1}{ \nu (U_n)} \int_{U_n} f d(\nu)-WMV^U_\mu(f)\right\} = 0$$
since $\frac{t\nu(U_n) }{ (t\nu + (1-t)\nu')(U_n)} \in [0;1]$ and $\lim_{n\rightarrow +\infty}\frac{1}{ \nu (U_n)} \int_{U_n} f d(\nu)=WMV_\nu^U(f)=WMV_\mu^U(f),$
and finally that $\lim_{n \rightarrow + \infty} \frac{t\nu(U_n) }{ (t\nu + (1-t)\nu')(U_n)} \left\{ \frac{1}{ \nu (U_n)} \int_{U_n} f d(\nu)-WMV^U_\mu(f)\right\} = 0$ the same way. 
Thus, $$WMV^U_{\nu''}(f)= \lim_{n \rightarrow +\infty}\frac{1 }{\nu''(U_n)} \int_{U_n} f d(\nu'')=WMV^U_\mu(f).$$
$\nu'' \in \mathcal{M}(\mu,U,F),$ thus $\mathcal{M}(\mu,U,F)$ is a convex cone.

\subsection{Asymptotic comparison of Radon measures}

We now turn to the number $\Theta$ that appeared in Proposition \ref{addition}.
In this section, $\mu$ and $\nu$ are fixed Radon measures and $U $ is a fixed sequence in $Ren_\mu \cap Ren_\nu.$ 
\begin{Proposition}
\begin{enumerate}
\item $ \Theta(\mu, \nu) \in [0;1].$
\item $\Theta (\mu, \nu) = 1 - \Theta(\nu, \mu).$
\end{enumerate}
\end{Proposition}

The proof is obvious.

\begin{Definition}\label{landau}
\begin{enumerate}
\item $\nu = o^U(\mu)$ if $\Theta(\mu,\nu)=1.$
\item $\nu = O^U(\mu)$ if $\Theta(\mu, \nu)>0.$
\item $\nu \sim^U \mu$ if $\Theta(\mu,\nu)=1/2.$
\end{enumerate}
\end{Definition}

Let us now compare three measures $\mu,$ $\nu$ and $\rho.$ 
The sequence $U$ is not precised now since it is a fixed arbitrary sequence.

\begin{Lemma}
Let $\mu$ and $\nu$ be two measures and let $U \in Ren_\mu\cap Ren_\nu.$
$$\Theta(\mu, \nu) = \frac{1 }{ 1 + \theta(\mu,\nu)}$$
where $\theta(\mu,\nu)= \lim_{n \rightarrow + \infty} \frac{\nu(U_n) }{ \mu(U_n)} \in \bar{\mathbb{R}}=[0;+\infty].$
\end{Lemma}

The proof is obvious.

\begin{Proposition}
Let $U \in Ren_\mu \cap Ren_\nu \cap Ren_\rho.$
\begin{enumerate}
\item $\theta(\mu, \rho) = \theta(\mu,\nu) \theta(\nu,\rho)$ if $(\theta(\mu,\nu), \theta(\nu,\rho))\notin \{(0;+\infty),(+\infty;0)\}$
\item $\Theta(\mu, \rho) = \frac{\Theta(\mu,\nu)\Theta(\nu,\rho) }{2 \Theta(\mu,\nu)\Theta(\nu,\rho) - \Theta(\mu,\nu) - \Theta(\nu,\rho)+1}$ if $(\Theta(\mu,\nu), \Theta(\nu,\rho))\notin \{(1;0),(0;1)\}$

\end{enumerate}
\end{Proposition}

\noindent
\textbf{Proof.}
Let $n\in \N.$ We have $\frac{\mu(U_n)}{(\mu + \nu)(U_n)} = \frac{1}{ 1+\frac{\nu(U_n)}{\mu(U_n)}},$ $\frac{\nu(U_n)}{ (\nu + \rho)(U_n)} = \frac{1}{ 1+\frac{\rho(U_n)}{\nu(U_n)}}$ and  $\frac{\mu(U_n)}{(\mu + \rho)(U_n)} = \frac{1 }{1+ \frac{\rho(U_n)}{\mu(U_n)}}.$ 
For the first part of the statement, $$\frac{\rho(U_n)}{\mu(U_n)} = \frac{\rho(U_n)}{\nu(U_n)}\frac{\nu(U_n)}{\mu(U_n)}.$$
(since these numbers are positive, the equality makes sense) Thus, if the limits are compatible, we get 1. taking the limits of both parts.
Then, we express each part  as : 
$$\frac{\rho(U_n)}{\mu(U_n)} = \frac{(\mu + \rho)(U_n) }{ \mu(U_n)}-1,$$
$$\frac{\rho(U_n)}{\nu(U_n)} = {(\nu + \rho)(U_n) }{\nu(U_n)}-1,$$
$$\frac{\nu(U_n)}{\mu(U_n)} = {(\mu + \nu)(U_n) }{ \mu(U_n)}-1,$$
and we get: 
$$\frac{\mu(U_n)}{(\mu + \rho)(U_n)}= \frac{\frac{\mu(U_n)}{(\mu+\nu)(U_n)} \frac{\nu(U_n)}{(\nu+\rho)(U_n)}}{2\frac{\mu(U_n)}{(\mu+\nu)(U_n)} \frac{\nu(U_n)}{(\nu+\rho)(U_n)}-\frac{\mu(U_n)}{(\mu+\nu)(U_n)}-\frac{\nu(U_n)}{(\nu+\rho)(U_n)} +1  }.$$
Taking the limit, we get 2.

\vskip 12pt

We recover by these results a straightforward extension of the comparison of the asymptotic behavour of functions. The notation chosen in \ref{landau} are chosen to show this correspondence. Through easy calculations of $\theta$ or $\Theta,$ one can easily see that, if $\mu, \nu$ and $\nu '$ are comparable measures, 
\begin{enumerate}
\item $(\mu \sim \nu) \wedge (\nu \sim \nu ') \Rightarrow (\mu \sim \nu ')$
\item $(\mu \sim \nu) \Leftrightarrow (\nu \sim \mu)$
\item $(\mu =o( \nu)) \Rightarrow (\mu = O(\nu))$
\item $(\mu =O( \nu)) \wedge (\nu =O(\nu ')) \Rightarrow (\mu =O(\nu '))$
\item $(\mu =o( \nu)) \wedge (\nu =o(\nu ')) \Rightarrow (\mu =o(\nu '))$
\item $(\mu =O( \nu ')) \wedge (\nu =O(\nu ')) \Rightarrow (\mu+\nu =O(\nu '))$
\end{enumerate}  
and other easy relations can be deduced in the same spirit. 

\subsection{Limits and mean value}

If $X$ is e.g. a connected locally compact, paracompact and not compact
manifold, equipped with a Radon measure $\mu$ such that $\mu(X)=+\infty,$
any exhaustive sequence $K=(K_{n})_{n\in\mathbb{N}}$ of compact subsets
of $X$ is such that $K\in Ren_{\mu}.$ In this setting, it is natural
to consider $\bar{X}=X\cup\infty$ the Alexandroff compactification
of $X.$

\begin{Theorem}\label{limit}
Let $f:X\rightarrow\mathbb{R}$ be a bounded measurable map which
extends to $\bar{f}:\bar{X}\rightarrow\mathbb{R},$ a continuous map
at $\infty.$ Then $WMV_{\mu}^{K}(f)=\bar{f}(\infty)$ for each exhaustive
sequence $K$ of compact subsets of $X.$

\end{Theorem}

\noindent
\textbf{Proof.}

We can assume that $\bar{f}(\infty)=0,$ in other words $$\lim_{x\rightarrow\infty}f(x)=0.$$The
sequence $(K_{n}^{c})_{n\in\mathbb{N}}$ gives a basis of neighborhood
of $\infty,$ thus $$\forall\epsilon'>0,\exists N'\in\mathbb{N},\forall n\geq N',sup_{x\in K_{n}^{c}}|f(x)|<\epsilon.$$
Moreover, since $\lim_{n\rightarrow+\infty}\mu(K_{n})=+\infty,$ $$\forall n_{0}\in\mathbb{N},\forall\epsilon''>0,\exists N''\in\mathbb{N},\forall n\geq N'',\mu(K_{n_{0}})<\epsilon''\mu(K_{n}).$$
Let $\epsilon>0.$ Let $\epsilon'=\frac{\epsilon}{2}.$ We set $n_{0}=N'$and
$\epsilon''=\frac{\epsilon}{2sup_{X}|f|}.$ Then, $\forall n\geq N=max(n_{0},N''),$
\begin{eqnarray*}|\int_{K_{n}}fd\mu| &\leq &\int_{K_{n}}|f|d\mu=\int_{K_{n_{0}}}|f|d\mu+\int_{K_{n}-K_{n_{0}}}|f|d\mu \\ & \leq &(sup_{X}|f|)\mu(K_{n_{0}})+\epsilon'\mu(K_{n}-K_{n_{0}}).\end{eqnarray*}
The second term is bounded by $\epsilon'\mu(K_{n})=\frac{\epsilon\mu(K_{n})}{2}$ and
we majorate the first term by $\epsilon''(sup_{X}|f|)\mu(K_{n})=\frac{\epsilon\mu(K_{n})}{2}.$ Thus
$$\forall\epsilon>0,\exists N>0,\forall n\geq N,|\frac{1}{\mu(K_{n})}\int_{K_{n}}fd\mu|\leq\epsilon,$$
and hence $WMV_{\mu}^{K}(f)=0.$

\vskip 12pt
As mentioned in introduction, we found no straightforward 
Beppo-Levy type theorem for mean values. The first counter-example we find is, 
for $X=\R$ and $\mu = \lambda$ the Lebesgue measure, an increasing sequence 
of $L^1(\lambda)$ which converges to ${1}_\R $ 
(uniformly on each compact subset of $\R ),$ e.g. 
the sequence $(e^{-\frac{x^2}{n}})_{n \in \N^*}.$ Let $K_n =[-n-1;n+1]$ and $K=(K_n)_{n\in \N}.$ We have $K\in Ren_\lambda,$
 $WMV_\lambda^K(1_\R)=1$ and $WMV_\lambda^K(e^{-\frac{x^2}{n}} )=0$ by 
Theorem \ref{limit}. 
We can only state the following theorem on uniform convergence:

\begin{Lemma}\label{borne}
Let $\mu$ be a measure on $X$ and let $U\in Ren_\mu .$ Let $f_1$ and $f_2$ be 
two functions in $\F_\mu^U(X,V)$ where $V$ in a sclctvs. 

Let $p$ be a 
norm on $V.$ If there exists $\epsilon \in \R_+^*$ such that $ \sup_{x \in X}\{p(f_1(x)-f_2(x))\} < \epsilon, $ then
$$ p(WMV_\mu^U(f_1))-\epsilon \leq p(WMV_\mu^U(f_2)) \leq p(WMV_\mu^U(f_1))+\epsilon.$$

\end{Lemma}

\vskip 12pt
\noindent
\textbf{Proof.}
Let $n \in \N.$
\begin{eqnarray*} p\left( \frac{1}{\mu(U_n)} \int_{U_n} f_2 d\mu \right)
& \leq &  \frac{1}{\mu(U_n)} \int_{U_n}p( f_1-f_2) d\mu + p\left( \frac{1}{\mu(U_n)} \int_{U_n} f_1 d\mu \right)\\
& \leq & \epsilon + p\left( \frac{1}{\mu(U_n)} \int_{U_n} f_1 d\mu \right) \end{eqnarray*}
We get the same way 
\begin{eqnarray*}  p\left( \frac{1}{\mu(U_n)} \int_{U_n} f_1 d\mu \right)-\epsilon& \leq &p\left( \frac{1}{\mu(U_n)} \int_{U_n} f_2 d\mu \right)\end{eqnarray*}
The result is obtained by taking the limit.

\begin{Theorem} \label{uniform}

Let $(f_n)_{n\in\N}\in \left(\F_\mu^U\right)^\N$ be a sequence which 
converges for uniform convergence on $X$ to a $\mu-$measurable map $f.$ Then 
\begin{enumerate}
\item $f\in \F_\mu^U .$
\item $WMV_\mu^U(f) = \lim_{n\rightarrow + \infty} WMV_\mu^U(f_n).$

\end{enumerate} 

\end{Theorem}

\noindent
\textbf{Proof.} Let $u_n = WMV_\mu^U(f_n).$ 

$\bullet$ \underline{Let us prove that} $(u_n)$ has a limit $u \in V.$

Let $p$ be a norm on $V.$ Let $\epsilon \in \R_+^*.$ There exists $N \in \N$ such that, for each $(n,m) \in \N^2,$ $$ sup_{x\in X} p(f_n - f_m) < \epsilon.$$
Thus, by Lemma \ref{borne} with $f_1 = 0$ and $f_2 = f_n - f_m,$   
$$ p(u_n -u_m) = p(WMV^U_\mu(f_n - f_m)) \leq \epsilon.$$
Thus, the sequence $(u_n)$ is a Cauchy sequence. 
Since $V$ is complete, the sequence $(u_n)$ has a limit 
$u \in V.$ 

$\bullet$ Moreover, we remember that  
$\forall \epsilon > 0, \forall (n,m) \in \N^2,$ 
\begin{eqnarray*} (sup_{x\in X} p(f_n - f) < \epsilon) \wedge (sup_{x\in X} p(f_m - f) < \epsilon) & \Rightarrow & sup_{x\in X} p(f_n - f_m) < 2\epsilon\\
& \Rightarrow & p(WMV^U_\mu (f_n - f_m)) < 2\epsilon \\
& \Rightarrow & p(u_n - u) < 2\epsilon \end{eqnarray*}

$\bullet$ \underline{Let us prove that} $u=\lim_{n \rightarrow +\infty }\frac{1}{\mu(U_n)}\int_{U_n}f d\mu .$
Let $(n,k) \in \N^2.$
\begin{eqnarray*}
p\left(\frac{1}{\mu(U_n)}\int_{U_n}f d\mu - u \right) & \leq & p\left(\frac{1}{\mu(U_n)}\int_{U_n}f d\mu - \frac{1}{\mu(U_n)}\int_{U_n}f_k d\mu \right) \\
&& + p\left( \frac{1}{\mu(U_n)}\int_{U_n}f_k d\mu - u_k \right) + p(u_k-u)
\end{eqnarray*}
Let $\epsilon \in \R_+^*.$ Let $K$ such that $\forall k > K,$ $sup_{x \in X}p(f-f_k)< \frac{\epsilon}{8}.$
Then $$p\left(\frac{1}{\mu(U_n)}\int_{U_n}f d\mu - \frac{1}{\mu(U_n)}\int_{U_n}f_k d\mu \right) < \frac{\epsilon}{8}$$
and $$p(u_k-u)<\frac{\epsilon}{4}.$$
Let $N$ such that for each $n > N,$ $$p\left( \frac{1}{\mu(U_n)}\int_{U_n}f_{K+1} d\mu - u_{K+1} \right)< \frac{\epsilon}{8}.$$
Then, by the same arguments, for each $k > K,$
$$p\left( \frac{1}{\mu(U_n)}\int_{U_n}f_k d\mu - u_k \right)< \frac{3\epsilon}{8}.$$
Gathering these inequalities, we get 
$$
p\left(\frac{1}{\mu(U_n)}\int_{U_n}f d\mu - u \right) < \frac{\epsilon}{8} +\frac{3\epsilon}{8}+ \frac{\epsilon}{4}=\epsilon.$$
This ends the proof of the theorem.

\subsection{Invariance of the mean value with respect to the Lebesgue measure}
In this section, $X = \R^m$ with $n \in \N^*,$ $\lambda$ is the Lebesgue measure, $K=(K_n)_{n \in \N}$ is the renormalization procedure defined by $$K_n = [-n-1;n+1]^n$$
and $L = (L_n)_{n \in \N}$ is the renormalization procedure defined by
$$ L_n = \{ x \in \R^n; ||x||\leq n+1 \}$$
where $||.||$ is the Euclidian norm. 
We note by $||.||_\infty$ the sup norm, and $d_\infty$ its associated distance. 
Let $v \in \R^n.$ We use the obvious notations $K+v=(K_n+v)_{n \in \N}$ and 
$L +v= (L_n+v)_{n \in \N}$ for the translated sequences. 
Let $(A,B)\in \mathcal{P}(X)^2.$
We note by  $A\Delta B = (A-B)\cup (B-A)$ the symmetric difference of subsets.

\begin{Proposition}\label{invariant0}
Let $v \in \R^m.$
Let $f \in \F_\lambda^K$ (resp. $f \in \F_\lambda^L$) be a bounded function. 
Let $U\in Ren_\lambda$ and $v\in \R^n.$ If  
$$\lim_{n \rightarrow +\infty}\frac{\lambda(U_n \Delta U_n+v)}{\lambda(U_n)} = 0,$$
\begin{enumerate}

\item Then $f \in \F_\lambda^{U+v}$ (resp. $f \in \F_\lambda^{U+v}$) and
$WMV_\lambda^U(f) = WMV_\lambda^{U+v}(f).$

\item Let $f_v: x\mapsto f(x-v).$ Then $f \in \F_\lambda^U$  and $WMV_\lambda^U(f) = WMV_\lambda^{U}(f_v) .$
\end{enumerate}
\end{Proposition}
\noindent
\textbf{Proof.}
We first notice that the second item is a reformulation of the first item: 
by change of variables $x \mapsto x-v,$ $WMV_\lambda^{K+v}(f)=WMV_\lambda^{K}(f_v).$

Let us now prove the first item. 
Let $n \in \N.$ 
\begin{eqnarray*} &&\frac{1}{\lambda(U_n)}\int_{U_n} f d\lambda - \frac{1}{\lambda(U_n+v)}\int_{U_n+v} f d\lambda \\ & = & \frac{1}{\lambda(U_n)}\int_{U_n} f d\lambda - \frac{1}{\lambda(U_n)}\int_{U_n+v} f d\lambda \\
& = &\frac{1}{\lambda(U_n)}\left(\int_{U_n -(U_n+v)} f d\lambda - \int_{(U_n+v)-U_n} f d\lambda \right) \\
& = &\frac{1}{\lambda(K_n)}\left(\int_{U_n \Delta(U_n+v)} (1_{U_n -(U_n+v)}-1_{(U_n+v)-U_n})f d\lambda\right) \end{eqnarray*}
Let $M = sup_{\R^m}(|f|).$ Then
$$|\frac{1}{\lambda(U_n)}\int_{U_n} f d\lambda - \frac{1}{\lambda(U_n+v)}\int_{U_n+v} f d\lambda| \leq M\frac{\lambda(U_n \Delta U_n+v)}{\lambda(U_n)}.$$
Thus, we get the result.
\begin{Lemma}
$$\lim_{n \rightarrow +\infty}\frac{\lambda(K_n \Delta K_n+v)}{\lambda(K_n)} = 0$$
and 
$$\lim_{n \rightarrow +\infty}\frac{\lambda(L_n \Delta L_n+v)}{\lambda(L_n)} = 0$$
\end{Lemma}
\noindent
\textbf{Proof.}
We prove it for the sequence $K,$ and the proof is the same for the sequence $L.$
We have $K_n = (n+1) K_0$ thus $\lambda(K_n) = (n+1)^m\lambda(K_0)$ and 
$$\lambda(K_n \Delta K_n+v) = (n+1)^m \lambda(K_0 \Delta K_0 + \frac{1}{n+1} v).$$
Let $$A_n = \{x \in \R^m | d_\infty(x,\partial K_0) < \frac{2||v||_\infty}{n+1}.$$
We have $K_0 \Delta K_0 + \frac{1}{n+1} v \subset A_n$ and 
$\lim_{n \rightarrow +\infty} \lambda(A_n)=0.$ 
Thus, 
$$\lim_{n \rightarrow +\infty}\frac{\lambda(K_n \Delta K_n+v)}{\lambda(K_n)} = 0.$$ 

\vskip 12pt
\begin{Proposition}\label{invariant1}
Let $v \in \R^m.$
Let $f \in \F_\lambda^K$ (resp. $ f \in \F_\lambda^L$) be a bounded function. 

\begin{enumerate}

\item Then $f \in \F_\lambda^{K+v}$ (resp. $ f \in \F_\lambda^{L+v}$) and
$WMV_\lambda^K(f) = WMV_\lambda^{K+v}(f)$ (resp. $WMV_\lambda^L(f) = WMV_\lambda^{L+v}(f)$).

\item Let $f_v: x\mapsto f(x-v).$ Then $f \in \F_\lambda^K$ (resp. $ f \in \F_\lambda^L$) and $WMV_\lambda^K(f) = WMV_\lambda^{K}(f_v)$ (resp. $WMV_\lambda^L(f) = WMV_\lambda^{L}(f_v)$).
\end{enumerate}
\end{Proposition}
\noindent
\textbf{Proof.}
The proof for $K$ and $L$ is a straightforward application of Proposition \ref{invariant0} whch is valid thanks to the previous Lemma.
\vskip 12pt
Concerning mean values on $m-$dimensional vector spaces, we must remark that 
the difference between two finite weak mean values of a same function 
$f$ can be huge. 
For $m=2,$ classical result of topology gives: 
$$ L_n \subset K_n \subset L_{E(\sqrt{2}n + \sqrt{2})}.$$
Let $f\in \F_\lambda^{K} \cap \F_\lambda^{L}$ be a positive function. 
For $n \in \N,$
\begin{eqnarray*}
\frac{1}{\lambda(K_n)}\int_{K_n}fd\lambda - \frac{1}{\lambda(L_n)}\int_{L_n}fd\lambda & = &
\frac{1}{\lambda(K_n)}\left\{\int_{K_n}fd\lambda - \frac{(n+1)^2}{\pi(n+1)^2}\int_{K_n}1_{L_n}fd\lambda \right\} \\
& = & \frac{1}{\lambda(K_n)}\left\{\int_{K_n}(1 - \frac{1}{\pi}1_{L_n})fd\lambda \right\}
\end{eqnarray*}
This shows that there can be a difference between $WMV^K_\lambda$ and $WMV^L_\lambda.$
\subsection{Example: the mean value induced by a smooth Morse function} \label{smorse}
In this example, $X$ is a smooth, locally compact, 
paracompact, connected, oriented and non compact manifold of dimension $n\geq 1$ equipped 
with a measure $\mu$ induced by a volume form $\omega$ and 
a Morse function $F : X \rightarrow \R$ such that $$\forall a \in \R, \quad \mu(F \leq a) < +\infty.$$
For the theory of Morse functions we refer to \cite{Mil}.
Notice that there exists some value $A$ such that $\mu(F<A) > 0.$ Notice that we can have $\mu(X)\in ]O;+\infty].$

\begin{Definition}
Let $f : X \rightarrow V$ be a smooth function into a sclctvs $V.$ Let $t \in [A;+\infty[. $ We define 
$$ I^f_\mu(f,t) = \frac{1}{\mu(F \leq t)} \int_{\{F \leq t\}} f(x)d\mu(x)$$
and, if the limit exists, 
$$WMV^F_\mu(f)= \lim_{t\rightarrow +\infty} I^F_\mu(f,t).$$

\end{Definition}

Of course this definition is the ``continuum'' version of the ``sequential'' definition \ref{def1}. 
If $V$ is metrizable, for any increasing sequence 
$(\alpha_n)_{n \in \N}\in [A;+\infty[^\N$ 
such that $\lim_{n \rightarrow +\infty}\mu\{F\leq \alpha_n\}\geq \mu(X),$ 
setting $U_n =\{F\leq \alpha_n\}$, $$WMV^U_\mu(f)=WMV^F_\mu(f)$$
and conversely $WMF^F_\mu(f)$ exists if $WMF^U_\mu(f)$ exists and does not depend on the choice of the 
sequence $(\alpha_n)_{n\in\N}.$ 

\vskip 12pt 
Moreover, since $F$ is a Morse function, it has isolated critical points and changing $X$ into $X-C,$ 
where $C$ is the set of critical points of $F,$ for each $t \in [A;+\infty[,$ $$\{F=t\}=F^{-1}(t)$$
is a $(n-1)-$dimensional manifold (disconnected or not). 
The first examples that we can give are definite positive 
quadratic forms on a vector space in which $X$ is embedded.

\subsection{Application: homology as a mean value} \label{Hodge}

Let $M$ be a finite dimensional manifold quipped with a 
Riemannian metric $g$ and the corresponding Laplace-Beltrami operator $\Delta,$
and with finite dimensional de Rham cohomology space $H^*(M,R).$
One of the standard results of Hodge theory is the onto and one-to-one map 
between 
$H^*(M,R)$ and the space of $L^2-$harmonic forms $\mathcal{H} $ made by integration over simplexes:
\begin{eqnarray*} I :& \mathcal{H} \rightarrow & H^*(M,R)  \\
			& \alpha \mapsto & I(\alpha)
\end{eqnarray*}
where $$I(\alpha): s \hbox{ simplex } \mapsto I(\alpha)(s)=\int_{s}\alpha.$$ 
We have  assumed here that the order of the simplex was the same as the order of the harmonic form.
This is mathematically coherent stating $\int_s\alpha = 0$ if $s$ and $\alpha$ do not gave the same order. 
Let $\lambda$ be the Lebesgue measure on $\mathcal{H}$ with respect to the scalar product induced by 
the $L^2-$scalar product. Let $U=(U_n)_{n\in \N}$ be the sequence 
of Euclidian balls centered at $0$ such that, for each $n\in \N,$ the ball $U_n$ is of radius $n.$
\begin{Proposition}
Assume that $H^*(M,R)$ is finite dimensional
Let $s$ be a simplex. Let $$\varphi_s = \frac{|I(.)(s)|}{1+|I(.)(s)|}.$$
The cohomology class of $s$ is null if and only if  $$WMV^U_\lambda(\varphi_s) = 0.$$
\end{Proposition}
\vskip 12pt
\noindent
\textbf{Proof.}

$\bullet$ If the cohomology class of $s$ is null, $\forall \alpha \in \mathcal{H},$ $\int_s \alpha = 0$ 
thus $\varphi_s(\alpha)=0.$
Finally, $$WMV^U_\lambda(\varphi_s) = 0.$$

$\bullet$ If the cohomology class of $s$ is not null, 
let $\alpha_s$ be the corresponding element in $\mathcal{H}.$ We have $\int_s\alpha_s=1.$
Let $\pi_s $ be the projection onto the 1-dimensional vector space spanned by $\alpha_s.$
Let $n \in \N^*.$ Let 
$$V_n = \left\{ \alpha \in U_n \hbox{ such that } |\int_s\alpha|>\frac{1}{2}.\right\}.$$
Then, $$V_n = U_n \cap \pi_s^{-1}\left([-1;1].\alpha_s\right).$$ 
Moreover, $$\inf_{\alpha \in V_n} \varphi_s(\alpha) = \frac{1}{2}$$
and $$\lambda(V_n) > \lambda(U_{n-1}).$$
Then \begin{eqnarray*} \int_{U_n}\varphi_s d\lambda & \geq & \int_{V_n}\varphi_s d\lambda\\
& \geq & \frac{\lambda(V_n)}{2}\end{eqnarray*}
Thus
\begin{eqnarray*}
WMV_\lambda^U(\varphi_s) & = & 
\lim_{n \rightarrow +\infty} \frac{1}{\lambda(U_n)} \int_{U_n}\varphi_s d\lambda \\
& \geq & \lim_{n \rightarrow +\infty}\frac{\lambda(V_n)}{2\lambda(U_n)}\\
& \geq & \lim_{n \rightarrow +\infty}\frac{\lambda(U_{n-1})}{2\lambda(U_n)}=\frac{1}{2}\\
& \neq & 0
\end{eqnarray*}

  \section{The mean value on infinite products} \label{IP}

\subsection{Mean value on an infinite product of measured spaces} \label{sprod}

Let $\Lambda$ be an infinite (countable, continuous or other) set of indexes. Let $(X_\lambda,\mu_\lambda)_{\lambda \in \Lambda}$ or for short $(X_\lambda)_\Lambda$ be a family of measured spaces as before. We assume that, on each space $X_\lambda,$ we have \textbf{fixed} a sequence $U_\lambda \in Ren_{\mu_\lambda}.$ Let $X_c=\prod_{\lambda\in\Lambda}X_\lambda$ be the cartesian product of the sequence $(X_\lambda)_\Lambda .$ 
\begin{Definition}
Let $f \in C^0(X)$ for the product topology. $f$ is called \textbf{cylindrical} if and only if there exists $\tilde{\Lambda}$ a finite subset of $\Lambda$ and a map $\tilde{f} \in C^0(\prod_{\lambda \in \tilde{\Lambda}}X_\lambda)$ such that $$\forall (x_\lambda)_{\Lambda}\in X,  f((x_\lambda)_\Lambda)=\tilde{f}((x_\lambda)_{\tilde{\Lambda}}).$$  
Then, we set, if $\tilde{f}\in \F^{\prod_{\lambda \in \tilde{\Lambda}}U_\lambda}_{\bigotimes_{\lambda \in \tilde{\Lambda}}\mu_\lambda},$ $$WMV(f) = WMV^{\prod_{\lambda \in \tilde{\Lambda}}U_\lambda}_{\bigotimes_{\lambda \in \tilde{\Lambda}}\mu_\lambda}(\tilde{f}).$$
We set the notation : $f \in \F.$
(here, subsidiary notations are omitted since the sequence of 
measures and the sequences of renormalization are fixed in this section)

We set the notation : $f \in \F.$ 
\end{Definition}

Notice that if we have $\tilde{\Lambda} \subset \tilde{\Lambda_0}$ with the notations used in the definition, 
since $f$ is constant with respect to the variables $x_\lambda$ indexed 
by $\lambda \in \tilde{\Lambda_0} - \tilde{\Lambda},$ 
the definition of $WMV(f)$ does not depend on the choice of $\tilde{\Lambda},$ which makes it coherent.  
\begin{Theorem}
Let $f$ be a cylindrical function associated to the finite set of indexes $\tilde{\Lambda}=\{\lambda_1,...,\lambda_n\}$ and to the function $\tilde{f}\in C^0(\prod_{\lambda \in \tilde{\Lambda}}X_\lambda).$ 
\begin{enumerate}
\item Let $\lambda \in \tilde{\Lambda}$. Let us fix $U_\lambda \in Ren_{\mu_\lambda}.$ Then $\prod_{\lambda \in \tilde{\Lambda}}U_\lambda = \left(\prod_{\lambda \in \tilde{\Lambda}}(U_\lambda)_n\right)_{n\in\N} \in Ren_{\bigotimes_{\lambda \in \tilde{\Lambda}}\mu_\lambda}.$

\item If both sides are defined, for each scalar-valued map $f=f_{\lambda_1}\otimes ... \otimes f_{\lambda_n} \in \F_{\mu_{\lambda_1}}^{U_{\lambda_1}}\otimes ... \otimes \F_{\mu_{\lambda_1}}^{U_{\lambda_n}} ,$
$$ WMV^{\prod_{\lambda \in \tilde{\Lambda}}U_\lambda}_{\bigotimes_{\lambda \in \tilde{\Lambda}}\mu_\lambda}(\tilde{f})=\prod_{\lambda \in \tilde{\Lambda}} WMV_{\mu_\lambda}^{U_\lambda}(f_\lambda).$$

\end{enumerate}
\end{Theorem}

For convenience of notations, we shall write $WMV(\tilde{f})$ instead of $WMV^{\prod_{\lambda \in \tilde{\Lambda}}U_\lambda}_{\bigotimes_{\lambda \in \tilde{\Lambda}}\mu_\lambda}(\tilde{f}).$
Let us now consider an arbitrary map $f:X\rightarrow V$ which is not cylindrical ($V$ is a sclctvs). 
Theorem  \ref{uniform} gives us a way to extend the notion of mean value by uniform convergence 
of sequences of cylindrical maps. 
But we shall not only do this for $X,$ but for classes of functions 
defined on a class of subset of $X.$ These classes are the following ones

\begin{Definition} \label{admissible}
Let $\mathcal{D}\subset X.$ The domain $\mathcal{D}$ is called \textbf{admissible} if and only if $$\forall x \in \mathcal{D}, \forall \tilde{\Lambda} \hbox{ finite subset of }\Lambda, \forall n\in \N, \left(\bigotimes_{\lambda \in \tilde{\Lambda}}\mu_\lambda\right) \left( \left(\prod_{\lambda \in \tilde{\Lambda}}U_{\lambda,n} \right)- \mathcal{D}_{\tilde{\Lambda},n,x}\right)= 0,$$
where $$\mathcal{D}_{\tilde{\Lambda},n,x}=\left\{u \in \prod_{\lambda \in \tilde{\Lambda}}U_{\lambda,n}|\exists x'\in\mathcal{D}, (\forall \lambda \in \tilde{\Lambda}, x'_\lambda = u_\lambda)\wedge(\forall \lambda \in \Lambda-\tilde{\Lambda}, x'_\lambda = x_\lambda)\right\}.$$
\end{Definition}

\begin{Definition}
Let $\mathcal{D}$ be an admissible domain. 
A function $f:\mathcal{D} \rightarrow V$ is \textbf{cylindrical} if its 
value depends only on a finite number of coordinates indexed by a fixed finite subset
of $\Lambda.$
\end{Definition}

The mean value of a cylindrical function $f$ comes immediately, since its trace defined 
on $\prod_{\lambda \in \tilde{\Lambda}}U_{\lambda,n}$ up to a subset of measure 0.
\begin{Theorem} \label{uniform2}
Let $V$ be a sclctvs. Let $f : \mathcal{D} \rightarrow V$ be the uniform limit of a 
sequence $(f_n)_{n \in \N}$ of cylindrical functions on $\mathcal{D}$ 
with a mean value on $\mathcal{D}.$
Then, 

\begin{enumerate}
\item the sequence $(WMV(f_n))_{n \in \N}$ has a limit.
\item This limit does not depend on the sequence $(f_n)_{n\in \N}$ but only on $f.$
\end{enumerate}
\end{Theorem}

\vskip 12pt
\noindent
\textbf{Proof.}

\noindent
Let $u_n = WMV_\mu^U(f_n).$ 

$\bullet$ \underline{Let us prove that} $(u_n)$ has a limit $u \in V.$

Let $p$ be a norm on $V.$ Let $\epsilon \in \R_+^*.$ There exists $N \in \N$ such that, for each $(n,m) \in \N^2,$ $$ sup_{x\in X} p(f_n - f_m) < \epsilon.$$
Thus, by Lemma \ref{borne} with $f_1 = 0$ and $f_2 = f_n - f_m,$   
$$ p(u_n -u_m) = p(WMV^U_\mu(f_n - f_m)) \leq \epsilon.$$
Thus, the sequence $(u_n)$ is a Cauchy sequence. 
Since $V$ is complete, the sequence $(u_n)$ has a limit 
$u \in V.$

$\bullet$ Now, let us consider another sequence $(f'_n)$ of cylindrical functions which converge uniformly to $f.$ 
In order to finish the proof of the theorem,  

\centerline{\underline{let us prove that} $u = \lim_{n \rightarrow +\infty} WMV^U_\mu(f'_n).$}

Let $n \in \N.$ We define $$f''_n =\left\{\begin{array}{l} f'_{\frac{n}{2}}  \hbox{ if n is even} \\
f_{\frac{n-1}{2}}  \hbox{ if n is odd}
\end{array} \right. .$$ 
This sequence again converges uniformly to $f$, and is hence a Cauchy sequence. By the way, the sequence $(WMV_\mu^U(f''_n)_{n \in \N}$ has a limit $u' \in V.$ Extracting the sequences $(f_n)_{n\in \N} = (f''_{2n+1})_{n\in \N} $ and $(f'_n)_{n\in \N} = (f''_{2n+1})_{n\in \N} $ we get
$$ u' = \lim_{n\rightarrow +\infty} WMV_\mu^U(f''_n) = \lim_{n \rightarrow +\infty} WMV_\mu^U(f_n) = u$$
and 
$$ \lim_{n \rightarrow +\infty} WMV_\mu^U(f'_n) = \lim_{n \rightarrow +\infty} WMV_\mu^U(f''_n)= u.$$

\vskip 12pt
By the way, the following definition is justified:
\begin{Definition}
Let $V$ be a sclctvs. Let $f : \mathcal{D} \rightarrow V$ be the uniform limit of a 
sequence $(f_n)_{n \in \N}$ of cylindrical functions on $\mathcal{D}$ 
with a mean value on $\mathcal{D}.$
Then, $$ WMV_\mu^U(f) = \lim_{n \rightarrow +\infty} WMV_\mu^U (f_n).$$
\end{Definition}

\noindent
Trivially, the map $WMV_\mu^U$ is linear 
as well as in the context of Proposition \ref{linear}.
\subsection{Application: the mean value on marked infinite configurations}

Let $X$ be a locally compact and paracompact manifold, 
orientable, and let $\mu$ be a measure on $X$ 
induced by a volume form.
In the following, we have either

- if $X$ is compact, setting $x_0 \in X,$ 
$$ \Gamma = \{ (u_n)_{n\in \N} \in X^\N | \lim u_n = x_0 \hbox{ and } \forall (n,m)\in \N^2, n\neq m \Rightarrow u_n \neq u_m\} $$ 

- if $X$ is not compact, setting $(K_n)_{n\in\N}$ an exhaustive sequence of compact subspaces of $X,$ 

$$O\Gamma = \{ (u_n)_{n\in \N} \in X^\N |  \forall p \in \N,  |\{u_n;n\in\N\} \cap K_p|<+\infty \hbox{ and } \quad \quad \quad \quad \quad \quad$$ 
$$ \quad \quad \quad \quad \quad \quad \quad \quad \quad \forall (n,m)\in \N^2, n\neq m \Rightarrow u_n \neq u_m\} $$
The first setting was first defined by Ismaginov, Vershik, Gel'fand and Graev,
see e.g. \cite{Ism} for a recent reference, and the 
second one has been extensively studied 
by Albeverio, Daletskii, Kondratiev, Lytvynov, see e.g. \cite{A}.
Alternatively, $\Gamma$ can be seen as a set of countable sums of Dirac measures
equipped with the topology of vague convergence.

For the following, we also need the set of ordered finite $k-$configurations: 
$$O \Gamma^k = \{ (u_1,...,u_k) \in X^k | \forall (n,m)\in \N^2, (1\leq n < m \leq k) \Rightarrow (u_n \neq u_m) \} $$
Assume now that $X$ is equipped with a Radon measure $\mu.$
One can notice that given $x \in \Gamma$ and a cylindrical function $f,$

Let us \textbf{fix} $U \in Ren_\mu .$ Notice first that for each $(n,k) \in \N^*\times \N ,$ 
$$ \mu^{\otimes n}(U_k^n) = \mu^{\otimes n}\left(
\{(x_1,...,x_n)\in U_k^n | \forall (i,j), (0\leq i < j \leq n)\Rightarrow (x_i \neq x_j) \} \right) .$$
In other words, the set of $n-$uples for which there exists two coordinates that are equal is of measure $0.$ This shows that $O\Gamma$ is an admissible domain in $X^\N,$ and enables us to write, for a bounded cylindrical function $f,$  

$$ WMV_\mu^U(f)=WMV_{\mu^{\otimes Ord(f)}}^{U^{Ord(f)}}(\tilde f)$$
since $\tilde f$ is defined up to a subset of measure $0$ on each $U_k^{Ord(f)},$ for $k \in \N.$
By the way, Theorem \ref{uniform2} applies in this setting 
Notice also that we the normalization sequence $U$ on $O\Gamma$ is induced from 
the normalization sequence on $X^\N .$ This implies heuristically that 
cylindrical functions with a weak mean value with respect to $U$ are somewhat 
small perturbations of functions on $X^\N.$
This is why we can modify the sequence $U$ on $O\Gamma$ the following way:
let $\varphi: \R_+ \rightarrow \R_+^*$ be a function such that 
$\lim_{x \rightarrow +\infty}\varphi =0.$ 
Then, if $f$ is a cylindrical function on $O\Gamma,$ we set 
$$U^n_\varphi = U^n - \{(x_i)_{1 \leq i \leq n} | \exists (i,j) \hbox{ such that } i < j \wedge d(x_i,x_j)<\varphi(n)\}.$$

\section{Mean value for heuristic Lebesgue measures} 

Any Fr\'echet space is the projective limit of a sequence of Banach spaces. 
Thus, any Fr\' echet space can be embedded in a Banach space $B,$ 
with continuous inclusion and density. We choose here 
to replace the Banach space $B$ by a Hilbert space $H$ in order to get (orthogonal)
canonical complementary subspaces.

\begin{Definition}
A \textbf{normalized Fr\'echet space}  is a pair $(F,H),$ where 

\begin{enumerate}
\item $F$ is a Fr\'echet space,

\item $H$ is a Hilbert space, 

\item $F \subset H$ and

\item $F$ is dense in $H.$
\end{enumerate}

\end{Definition}

Another way to understand this definition is the following: 
we choose a pre-Hilbert norm on the Fr\'echet space $F.$ 
Then, $H$ is the completion of $F.$

\begin{Definition}\label{dh}
Let $V$ be a sclctvs. A function $f : F \rightarrow V$ is \textbf{cylindrical} 
if there exists $F_f , $  a finite dimensional affine subspace of $F,$ for which, if $\pi$ is the orthogonal projection, 
$\pi : F \rightarrow F_f$ such that $$\forall x \in F, f(x) = f\circ \pi(x).$$  

\end{Definition}

\begin{Proposition}
Let $(f_n)_{n \in \N}$ be a sequence of cylindrical functions. There exists an
unique sequence $(F_{f_n})_{n\in \N}$ creasing for $\subset,$ for which 
$\forall m \in \N,$, $F_{f_m}$ is the minimal affine space for which $$\forall n \leq m, f_n \circ \pi_m = f_n.$$ 
\end{Proposition}

\noindent
\textbf{Proof.}
We build it by induction: 

$\bullet$ $F_{f_0}$ is the minimal affine subspace of $F$ for which
Definition \ref{dh} applies to $f_0.$

$\bullet$ Let $n \in \N.$ Assume that we have constructed $F_{f_n}.$ Let $\tilde{F}$
be the minimal affine subspace of $F$ for which
Definition \ref{dh} applies to $f_{n+1}.$ 
We set $$F_{f_{n+1}} = F_{f_n} + \tilde{F}.$$
(recall that it is the minimal affine subspace of $F$ which contains both $F_{f_n}$ and $\tilde{F}.)$
If $\tilde\pi$  and $\pi_{n+1}$ are the orthogonal projections into 
$\tilde F$ and $F_{f_{n+1}},$ that 
$$f_{n+1} = f_{n+1}\circ\tilde\pi = f_{n+1}\circ \pi_{n+1}.$$   
This ends the proof. $\qed$

We now develop renormalization procedures on $F$ inspired from 
section \ref{sprod}, using orthogonal projections to 1-dimensional vector subspaces. 
In these approaches, a finite dimensional Euclidian space is 
equipped with its Lebesgue measure noted by $\lambda$ in any dimension. 
The Euclidian norms are induced by the pre-Hilbert norm on $F$ for any 
finite dimensional vector subspace of $F.$  

\subsection{Mean value by infinite product}

Let $f$ be a bounded function which is the uniform limit of a sequence of 
cylindrical functions $(f_n)_{n \in \N}.$  
Here, an orthonormal basis 
$(e_k)_{k \in \N}$ is obtained by induction, completing at each 
step an orthonormal basis of $F_{f_n}$ by an orthonormal basis 
of $F_{f_{n+1}}.$ Thus we can identify $F$ with a subset $\mathcal{D}$ of $\R^\N$
which is invariant under change of a finite number of coordinates. 
This qualifies it as admissible since, with the notations used in Definition \ref{admissible}, 
    $$\left(\prod_{\lambda \in \tilde{\Lambda}}U_{\lambda,n} \right)- 
\mathcal{D}_{\tilde{\Lambda},n,x} = \emptyset,$$
for any set of renormalization procedures in $\R^\N$ as defined in section \ref{sprod}.  
So that, Theorem \ref{uniform2} applies. 
We note by $$ WMV_\lambda(f)$$ this value.
We remark that we already know by Theorem \ref{uniform2} that this mean 
value does not depend on the sequence $(f_n)_{n \in \N}$ 
\textbf{only once the sequence} $(F_{f_n})_{n \in \N}$ \textbf{is fixed}. In other words, 
two sequences $(f_n)_{n \in \N}$ and $(f'_n)_{n \in \N}$ which converge 
uniformly to $f$ a priori lead to the same mean value if $F_{f_n}=F_{f'_n}$ 
(maybe up to re-indexation). From heuristic calculations, it seems to come from 
the choice of the renormalization procedure, 
which is dependent on the basis chosen, 
more than from the sequence $(F_{f_n})_{n \in \N}.$ 
The problem would be solved if we did not get technical difficulties 
to replace the cubes $[-n-1;n+1]^k$, 
for $(n,k)\in \N\times\N^*,$ by an Euclidian ball. 
Further investigations are in progress.

\subsection{Invariance}

We notice three types of invariance: scale invariance,
translation invariance 
and invariance under the orthogonal (or unitary) group. 

\begin{Proposition}
Let $\alpha \in \N^*.$ Let $f$ be a function on $F$ with mean value. 
Let $f_\alpha: x\in F \mapsto f(\alpha x).$ Then $f_\alpha$ has a 
mean value and $$WMV_\lambda(f_\alpha) = WMV_\lambda(f).$$
\end{Proposition}  

\noindent
\textbf{Proof.}

Let $(f_n)_{n\in \N}$ be a sequence which converges uniformly to $f.$
Then, with the notations above, the sequence $((f_n)_\alpha)_{n\in \N}$
converges uniformly to $f_\alpha.$
Let $m = dim(F_{f_n}).$ 
$$WMV_\lambda((f_n)_\alpha) = WMV_{\alpha^{-m} \lambda}(f_n) = WMV_\lambda(f_n).$$
by proposition \ref{scaling} and remarking that 
for the fixed renormalization sequence above, 
this change of variables consists in extracting a subsequence of renormalization. 
Thus, taking the limit, we get $$WMV_\lambda(f_\alpha) = WMV_\lambda(f).$$

\begin{Proposition}
Let $v \in F.$ Let $f$ be a function on $F$ with mean value. 
Let $f_v: x\in F \mapsto f(x+v).$ 
Then $f_v$ has a mean value and $$WMV_\lambda(f_v)=WMV_\lambda(f).$$ 
\end{Proposition}

\noindent
\textbf{Proof.}
Let $(f_n)_{n\in \N}$ be a sequence which converges uniformly to $f.$
Let $v_n = \pi_n(v) \in F_{f_n}.$ We have $(f_n)_v=(f_n)_{v_n}.$
Then, 
\begin{eqnarray*}
WMV_\lambda(f_v) & = & \lim_{n \rightarrow +\infty} WMV_\lambda ((f_n)_{v_n}) \\
& = & \lim_{n \rightarrow +\infty} WMV_\lambda (f_n) \hbox{ by Proposition \ref{invariant1}}\\
& = & WMV_\lambda(f)
\end{eqnarray*}

\begin{Proposition}
Let $U_F$ be the group of unitary operators of $H$ 
which restricts to a bounded map $F \rightarrow F$ 
and which inverse restricts also. Let $u \in U_F.$ Let $f$
be a map with mean value. Then $f \circ u$ has a mean value and 
$$ WMV_\lambda(f \circ u)=WMV_\lambda(f).$$ 
\end{Proposition}

This last proposition becomes obvious after remarking that 
we transform the sequence $(F_{f_n})_{n\in \N}$ into the orthogonal sequence
$$ \left(u^{-1}(F_{f_n})\right)_{n\in \N} = (F_{f_n \circ u})_{n \in \N}.$$
This remark shows that we get the same mean value for $f\circ u$ as for $f$ by changing the orthogonal sequence. 

\subsection{Final remark: Invariance by restriction}

Let $G$ be a vector subspace of $F$ such that $$\bigoplus_{n\in \N} F_{f_n} \subset G.$$
As a consequence, if $g$ is the restriction of $f$ to $F_1,$ the sequence 
$(f_n)_{n\in \N}$ of cylindrical functions on $F$ restricts to a sequence 
$(g_n)_{n\in \N}$ of cylindrical functions on $G.$ Then, for uniform convergence, 
$$\lim_{n\rightarrow +\infty} g_n = g$$
and for fixed $n \in \N$ we get through restriction to $F_{f_n},$
$$WMV_\lambda(g_n) = WMV_\lambda(f_n).$$
Taking the limit, we get   $$WMV_\lambda(g) = WMV_\lambda(f).$$
This shows the restriction property announced in the introduction.

\end{document}